\documentclass[12pt]{amsart}
\usepackage{graphicx}
\usepackage{amsfonts}
\usepackage{amssymb}
\usepackage{latexsym}
\usepackage{amscd}

\renewcommand{\marginpar}[1]{}
\catcode`\@=12

\def\Empty{}
\newcommand\oplabel[1]{
  \def\OpArg{#1} \ifx \OpArg\Empty {} \else
  	\label{#1}
  \fi}
		
\long\def\realfig#1#2#3#4{
\begin{figure}[htbp]
\centerline{\psfig{figure=#2,width=#4}}
\caption[#1]{#3}
\oplabel{#1}
\end{figure}}

\newcommand{\comm}[1]{}

\input{psfig}
\newtheorem{thm}{Theorem}[section]
\newtheorem{cor}[thm]{Corollary}
\newtheorem{lem}[thm]{Lemma}
\newtheorem{prop}[thm]{Proposition}

\newenvironment{pf}{\proof[\proofname]}{\endproof}
\newenvironment{pf*}[1]{\proof[#1]}{\endproof}
\usepackage{euscript}

\usepackage[OT2,OT1]{fontenc}

\newcommand{\cal}[1]{{\mathcal #1}}

\theoremstyle{definition}
\newtheorem{defn}{Definition}[section]

\theoremstyle{remark}
\newtheorem{rem}{Remark}[section]

\newcommand{\diam}{\operatorname{diam}}
\newcommand{\dist}{\operatorname{dist}}
\newcommand{\meas}{\operatorname{meas}}

\renewcommand{\mod}{\operatorname{mod}}
\newcommand{\tl}{\tilde}
\newcommand{\wtl}{\widetilde}
\newcommand{\eps}{\epsilon}

\newcommand{\aaa}[1]{{{\mathbf{#1}}}}

\newcommand{\hol}{{\aaa H}}

\renewcommand{\Im}{\operatorname{Im}}

\numberwithin{equation}{section}
\newcommand{\thmref}[1]{Theorem~\ref{#1}}
\newcommand{\propref}[1]{Proposition~\ref{#1}}

\newcommand{\lemref}[1]{Lemma~\ref{#1}}
\newcommand{\corref}[1]{Corollary~\ref{#1}}
\newcommand{\figref}[1]{Figure~\ref{#1}}

\newcommand{\cU}{{\cal U}}

\newcommand{\cW}{{\cal W}}

\newcommand{\cV}{{\cal V}}

\newcommand{\cP}{{\cal P}}

\newcommand{\cH}{{\cal H}}
\newcommand{\cR}{{\cal R}}
\newcommand{\cL}{{\cal L}}

\newcommand{\cE}{{\cal E}}

\newcommand{\CC}{{\Bbb C}}
\newcommand{\RR}{{\Bbb R}}
\newcommand{\TT}{{\Bbb T}}
\newcommand{\ZZ}{{\Bbb Z}}
\newcommand{\NN}{{\Bbb N}}
\newcommand{\DD}{{\Bbb D}}
\newcommand{\HH}{{\Bbb H}}
\newcommand{\QQ}{{\Bbb Q}}

\newcommand{\cren}{\cR_{\text cyl}}

\begin{document}
\addtolength{\evensidemargin}{-0.7in}
\addtolength{\oddsidemargin}{-0.7in}

\title[Siegel universality]
{Siegel disks and renormalization fixed points}
\author{Michael Yampolsky}
\begin{abstract}
In this note we construct hyperbolic fixed points for  cylinder renormalization of maps with 
Siegel disks.

\end{abstract}

\date{\today}
\thanks{This note was written in the Fall of 2003, and amended with a corollary of the new 
results of M.~Shishikura in the Summer of 2005.}

\maketitle

\section{Introduction}
A renormalization hyperbolicity conjecture has so far been established in two examples
of one-dimensional dynamical systems: the unimodal maps, in the works of Sullivan \cite{S1,Sul,MvS}, 
McMullen \cite{McM1,McM2}, and Lyubich \cite{Lyubich-feigenbaum,Lyubich-horseshoe};
and the critical circle maps, in the works of de~Faria and de~Melo \cite{dF1,dF2,dFdM1,dFdM2},
and the author \cite{Ya2,Ya3,Ya4}. In this paper we will add one more example to the list,
by constructing a hyperbolic fixed points of renormalization corresponding to Siegel disks.
Let us say that an irrational number $\theta$ is {\it golden} if it is represented by an 
infinite continued fraction 
 $$\theta=\cfrac{1}{N+\cfrac{1}{N+\cfrac{1}{\cdots}}},\text{ with }N\in\NN.$$
We introduce this notation by analogy with the golden mean $(\sqrt{5}-1)/2$ which is
expressed by such a fraction with $N=1$. 
It has long been known, that golden Siegel disks
in the quadratic family have self-similar scaling properties near the critical point,
explained by a renormalization hyperbolicity conjecture \cite{MN,Wi}.
MacKay and Persival \cite{MP} have conjectured in 1986, based on numerical evidence, the
existence of a hyperbolic renormalization horseshoe corresponding to Siegel disks of 
analytic maps,
analogous to the Lanford's horseshoe for critical circle maps \cite{La1,La2}. 

In 1994 Stirnemann \cite{Stir} gave a computer-assisted proof of the existence of a
renormalization fixed point with a golden-mean Siegel disk.
In 1998, McMullen \cite{McM3} proved the asymptotic self-similarity of  golden
Sigel disks in the quadratic family. He constructed a version of renormalization 
based on holomorphic commuting pairs of de~Faria \cite{dF1,dF2} and
showed that the renormalizations of a quadratic polynomial with a golden Siegel disk
 near the critical point converge
to a fixed point geometrically fast. 
More generally, he constructed a renormalization horseshoe for bounded
type rotation numbers, and used renormalization to show that the Hausdorff dimension of
the corresponding quadratic Julia sets is strictly less than two.

In \cite{Ya3} we  introduced a new renormalization transformation $\cren$, which we called
the cylinder renormalization, and used it to prove the Lanford's Hyperbolicity
Conjecture for critical circle maps. The main advantage of $\cren$ over
the renormalization scheme based on commuting pairs is that this operator is analytic
in a Banach manifold of analytic maps of a subdomain of $\CC/\ZZ$. 
In this paper we  use McMullen's result 
to construct a fixed point $\hat f$ of $\cren$ with a golden-mean Siegel disk.
We further discuss the properties of $\cren$ at the fixed point $\hat f$, and, in particular, show
that its linearization is a compact operator having at least one eigenvalue outside the closed unit disk.

Finally, we use the new results of Inou and Shishikura \cite{Shi} to show that for sufficiently large values 
of $N$, this fixed point of $\cren$ is hyperbolic, and the dimension of its expanding subspace is
exactly one.

\begin{thm}
\label{main thm}
Let $$\theta_N=\cfrac{1}{N+\cfrac{1}{N+\cfrac{1}{\cdots}}},\text{ with }N\in\NN$$ be  a golden number. 
There exists a complex Banach space ${\aaa C}_{U(N)}$ whose elements (referred to as critical
cylinder maps) are analytic maps defined in a neighborhood
$U(N)$ of the origin, such that the following holds.
There exists a critical cylinder map
$\hat f_N\in{\aaa C}_{U(N)}$ with a Siegel disk $\Delta\Subset U(N)$ with rotation number $\theta_N$ for which:
\begin{itemize}
\item[(I)] the boundary of $\Delta$ is a quasicircle containing the critical point of $\hat f_N$;
\item[(II)] $\cren\hat f_N=\hat f_N$;
\item[(III)] the quadratic polynomial  $f(z)=e^{2\pi i\theta_N}z+z^2$
is infinitely cylinder renormalizable, and 
$$\cren^k f\to \hat f_N,$$
at a uniform geometric rate;
\item[(IV)] 
the cylinder renormalization $\cren$ is an analytic and compact operator mapping a neighborhood
of the fixed point $\hat f_N$ in ${\aaa C}_{U(N)}$ to ${\aaa C}_{U(N)}$. Its linearization at 
$\hat f_N$ is a compact operator, with at least one 
 eigenvalue with the absolute value greater than one;
\item[(V)] there exists a neighborhood $U_\infty$ of the origin which is contained in all $U(N)$ for 
large enough $N$.
\end{itemize}
\noindent
Moreover, there exists $N_0\in\NN$ such that for all $N\geq N_0$ we have:
\begin{itemize}
\item[(VI)] except for the one unstable eigenvalue, the spectrum of $\cren$ at $\hat f_N$ is compactly
contained inside the unit disk.
\end{itemize}

\end{thm}

\noindent
\begin{rem}
We note that the results (I)--(V) extend to show the existence of an invariant horseshoe
for renormalization of Siegel disks with rotation numbers $\theta$ of type bounded by any $B\in\NN$. The
result (VI) extends to the horseshoe of Siegel disks with bounded type for which every term
in the continued fraction expansion of $\theta_0$ is at least $N_0$. We restrict ourselves
to the case of a renormalization fixed point for simplicity of exposition.
\end{rem}

\subsection*{Acknowledgements} I would like to thank Xavier Buff for 
several useful discussions of the results of Inou and Shishikura, and
for suggesting a different argument for proving \thmref{thm cyl ren siegel}.
In addition, I would like to thank X.~Buff and A.~Ch{\'e}ritat for discussing
the Doaudy's Program of constructing positive measure Julia sets with me.

\section{Preliminary considerations.}
\label{section: preliminaries}

\noindent
{\bf Some notations.}
We use $\dist$ and $\diam$ to denote the Euclidean distance and diameter in $\Bbb C$.
We shall say that two real numbers $A$ and $B$ are $K$-commensurable for $K>1$ if
$K^{-1}|A|\leq |B|\leq K|A|$.
The notation $D_r(z)$ will stand for the Euclidean disk with the center at $z\in\Bbb C$ and
radius $r$. The unit disk $D_1(0)$ will be denoted $\Bbb D$. 
The plane $({\Bbb C}\setminus{\Bbb R})\cup J$
with the parts of the real axis not contained in the interval $J\subset \Bbb R$
removed  will be denoted ${\Bbb C}_J$.
By the circle $\TT$ we understand the affine manifold $\RR/\ZZ$, it is naturally identified
with the unit circle $S^1=\partial\DD$. The real translation $x\mapsto x+\theta$
projects to the rigid rotation by angle $\theta$, $R_\theta:\TT\to\TT$, the same
map on $S^1$ will be denoted $r_\theta$.
For two points $a$ and $b$ in the circle $\TT$
which are not diametrically opposite, $[a,b]$ will denote the shorter
of the two arcs connecting them. As usual, $|[a,b]|$ will denote the length of the 
arc. For two points $a,b\in \Bbb R$, $[a,b]$ will denote the closed interval with
endpoints $a$, $b$ without specifying their order.
The cylinder in this paper, unless otherwise specified will mean the affine manifold $\CC/\ZZ$.
Its equator is the circle $\{\Im z=0\}/\ZZ\subset \CC/\ZZ$.
A topological  annulus $A\subset \CC/\ZZ$ will be called an equatorial annulus, or an 
equatorial neighborhood, if it has a smooth boundary and contains the equator.

By ``smooth'' in this paper we will mean ``of class $C^\infty$'', unless another degree
of smoothness is specified. The notation``$C^\omega$'' will 
stand for ``real-analytic''.

We will sometimes use a symbol $\infty$ alongside the natural numbers, with the usual 
conventions $\infty>n$, $1/\infty=0$, and $\infty\pm n=\infty$ for $n\in\NN$. 

\medskip

\noindent
{\bf Renormalization of critical circle maps.} We are going to recall here very briefly the
way renormalization of critical circle maps is defined using commuting pairs \cite{FKS,ORSS}.
A detailed account of what follows may be found in \cite{Ya3}.
A critical circle mapping $f$ is a homeomorphism $\TT\to\TT$ of class $C^3$ with an only
critical point at $0$. The latter is further assumed to be non-flat, usually cubic.
Yoccoz \cite{Yoc} has shown that if such a mapping has an irrational rotation number
$\rho(f)$, then it is conjugate to the rigid rotation of the circle by angle $\rho(f)$
by a homeomorphic change of coordinate. Writing $\rho(f)$ as an infinite continued
fraction with positive terms
$$\rho(f)=\cfrac{1}{a_1+\cfrac{1}{a_2+\cfrac{1}{\cdots}}}$$
\noindent
(which we will further abbreviate as $[a_1,a_2,\ldots]$), we denote $\{p_n/q_n\}$ the sequence
of its convergents
$$p_n/q_n=[a_1,\ldots,a_n].$$
\noindent
As $p_n/q_n$ are best rational approximations of $\rho(f)$, the denominators $q_n$ are
closest return times of the critical point $0$: the arc $J_n=[0,f^{q_n}(0)]$ contains
no iterates $f^{j}(0)$ with $j<q_n$.
The first return map of $J_n\cup J_{n+1}$ is a piecewise defined mapping given by
$$R_{J_n\cup J_{n+1}}f=(f^{q_{n+1}}|_{J_n},f^{q_n}|_{J_{n+1}}).$$
This serves as a motivation for the following definition.

\begin{defn}
A  {\it  commuting pair} $\zeta=(\eta,\xi)$ consists of two 
$C^3$-smooth  orientation preserving interval homeomorphisms 
$\eta:I_\eta\to \eta(I_\eta),\;
\xi:I_{\xi}\to \xi(I_\xi)$, where
\begin{itemize}
\item[(I)]{$I_\eta=[0,\xi(0)],\; I_\xi=[\eta(0),0]$;}
\item[(II)]{Both $\eta$ and $\xi$ have homeomorphic extensions to interval
neighborhoods of their respective domains with the same degree of
smoothness, which commute, 
$\eta\circ\xi=\xi\circ\eta$;}
\item[(III)]{$\xi\circ\eta(0)\in I_\eta$;}
\item[(IV)]{$\eta'(x)\ne 0\ne \xi'(y) $, for all $x\in I_\eta\setminus\{0\}$,
 and all $y\in I_\xi\setminus\{0\}$.}
\end{itemize}
\end{defn}

\realfig{compair}{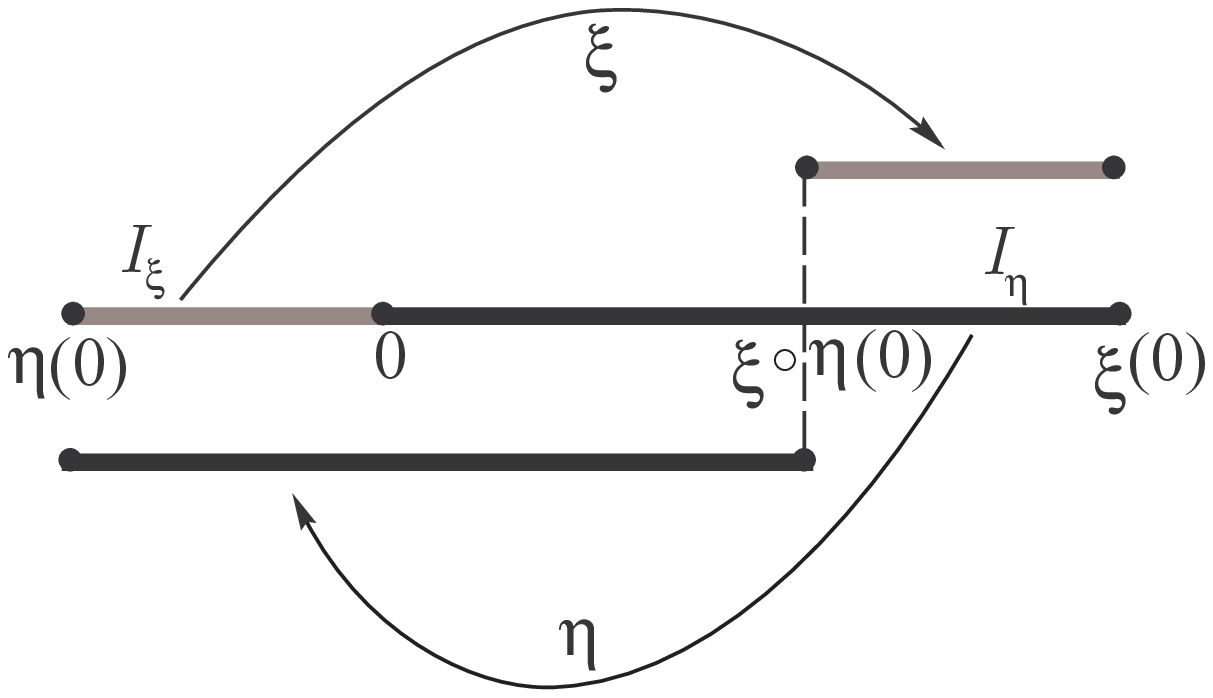}{A commuting pair}{8cm}

\noindent
The {\it height} $\chi(\zeta)$
of a critical commuting pair $\zeta=(\eta,\xi)$ is equal to $r$,
if 
$$0\in [\eta^r(\xi(0)),\eta^{r+1}(\xi(0))].$$
 If no such $r$ exists,
we set $\chi(\zeta)=\infty$, in this case the map $\eta|I_\eta$ has a 
fixed point.  For a pair $\zeta$ with $\chi(\zeta)=r<\infty$ one verifies directly that the
mappings $\eta|[0,\eta^r(\xi(0))]$ and $\eta^r\circ\xi|I_\xi$
again form a commuting pair.
For a commuting pair $\zeta=(\eta,\xi)$ we will denote by 
$\wtl\zeta$ the pair $(\wtl\eta|\wtl{I_\eta},\wtl\xi|\wtl{I_\xi})$
where tilde  means rescaling by the linear factor $\lambda=-{1\over |I_\eta|}$.

\begin{defn}
The {\it renormalization} of a real commuting pair $\zeta=(\eta,
\xi)$ is the commuting pair
\begin{center}
${\cal{R}}\zeta=(
\widetilde{\eta^r\circ\xi}|
 \widetilde{I_{\xi}},\; \widetilde\eta|\widetilde{[0,\eta^r(\xi(0))]}).$
\end{center}
\end{defn}

\noindent
For a pair $\zeta$ we define its {\it rotation number} $\rho(\zeta)\in[0,1]$ to be equal to the 
continued fraction $[r_0,r_1,\ldots]$ where $r_i=\chi({\cal R}^i\zeta)$. 
In this definition $1/\infty$ is understood as $0$, hence a rotation number is rational
if and only if only finitely many renormalizations of $\zeta$ are defined;
if $\chi(\zeta)=\infty$, $\rho(\zeta)=0$.

The non-rescaled pair $(\eta^r\circ\xi|I_\xi,\eta|[0,\eta^r(\xi(0))])$ will be referred to as the 
{\it pre-renormalization} $p{\cal R}\zeta$ of the commuting pair $\zeta=(\eta,\xi)$.
For a critical circle map $f$ as above, we set 
$$p\cR^nf=R_{J_n\cup J_{n+1}}f=(f^{q_{n+1}}|_{J_n},f^{q_n}|_{J_{n+1}}),\text{ and }\cR^nf=\widetilde{p\cR^nf}.$$

\noindent
A key object in the renormalization theory of commuting pairs developed by de~Faria
\cite{dF1,dF2} 
is the holomorphic commuting pair. This is an  analogue of the Douady-Hubbard's
polynomial-like map in the unimodal renormalization theory, which is defined as follows (cf. \figref{renormfig}):

\begin{defn}
An analytic commuting pair $\zeta=(\eta|_{I_\eta},\xi|_{I_\xi})$   extends to a {\it
holomorphic commuting
pair} $\cal H$ if there exist four simply-connected $\RR$-symmetric domains $V$, $U_1$, $U_2$, $U_3$ such that
\begin{itemize}
\item  $\bar U_1,\; \bar U_2,\; \bar U_3\subset V$,
 $\bar U_1\cap \bar U_2=\{ 0\}$; the sets
 $U_1\setminus U_3$,  $U_2\setminus U_3$, $U_3\setminus U_1$, and $U_3\setminus U_2$ 
 are nonempty, 
connected, and simply-connected, $U_1\cap \RR= I_\eta$, $U_2\cap \RR=I_{\xi}$;
\item mappings $\eta:U_1\to (V\setminus \RR)\cup\eta(I_\eta)$ and
 $\xi:U_2\to(V\setminus \RR)\cup\xi(I_\xi)$ are onto and
univalent;
\item $\nu\equiv \eta\circ\xi:U_3\to (V\setminus \RR)\cup{\nu(I_{U_3})}$ is a three-fold 
branched covering with a unique critical point at zero, where 
$I_{U_3}=U_3\cap \RR$.
\end{itemize}
\end{defn}

\noindent
One says that an analytic  commuting pair $(\eta,\xi)$ with
an irrational rotation number has
{\it complex {\rm a priori} bounds}, if all its renormalizations extend to 
holomorphic commuting pairs with {\it bounded moduli}:
$$\mod(V\setminus \overline{\cup U_i})>\mu>0.$$
Part of the significance of the complex {\it a priori} bounds is explained by the 
following:
\begin{prop}[\cite{Ya1}]\label{bounds compactness}
For $\mu\in(0,1)$ let $\hol(\mu)$ denote the space of holomorphic commuting pairs
${\cH}$, with $\mod (V\setminus\overline{\cup U_i})>\mu$,
$\min(|I_\eta|,|I_\xi|)>\mu$ and $\diam(V)<1/\mu$. Then the space $\hol(\mu)$ is
sequentially pre-compact with respect to the Carath{\'e}odory topology, with all the
limit points contained in $\hol(\mu/2)$.
\end{prop}
The existense of complex {\it a priori} bounds is a key analytic issue 
 of renormalization theory. In the case of critical circle maps it is settled
by the following theorem:

\begin{thm}
\label{complex bounds}There exist universal constants $\mu>0$ and $K>1$ such that
the following holds. 
Let $\zeta$ be an analytic critical commuting pair with an irrational rotation number. 
Then there exists $N=N(\zeta)$ such that for all
$n\geq N$ the  commuting pair $\cR^n\zeta$ extends to a holomorphic commuting pair
$\cH_n\in{\aaa H}(\mu)$.
Moreover, its range $V_n$ is a Euclidean disk, and the regions $U_i\cap(\pm \HH)$
are $K$-quasidisks.
\end{thm}

We first proved this theorem in \cite{Ya1} for commuting pairs 
$\zeta$ in the {\it Epstein class}. Our proof was later adapted by de~Faria and de~Melo \cite{dFdM2} 
to the case of a non-Epstein analytic commuting pair.

Finally, let us make the following note for future reference.
\begin{prop}
\label{dynamical partition}
If $f$ be a homeomorphism of a topological circle $S$ with $\rho(f)\in\RR\setminus\QQ$.
Fix a point $p\in S$ and let $A_n$ be the arc of the circle $S$ connecting $p$ to $f^{q_n}(p)$
and such that $f^{q_{n-1}}(p)\notin A_n$. Then for every $n\in\NN$
the topological arcs
\begin{equation}
\label{dyn part}
A_n, f(A_n),\ldots,f^{q_{n+1}-1}(A_n),A_{n+1},f(A_{n+1}),\ldots,f^{q_n-1}(A_{n+1})
\end{equation}
cover the whole circle and have disjoint interiors.
\end{prop}
\noindent
We will call the intervals (\ref{dyn part}) the {\it $n$-th dynamical partition} of $f$, given by the
orbit of $p$. 
\medskip


\noindent
{\bf Siegel quadratics.} Let $f(z)$ be a germ of an anlytic mapping at $0$, for which the
origin is a fixed point. We will concentrate on the case when $0$ is an irrationally
indifferent fixed point of $f$, that is, $f'(0)=e^{2\pi i\theta}$, 
with the rotation number 
$\theta\in\RR\setminus\QQ$. The germ $f$ is linearizable, if after a conformal change
of coordinate in a neighborhood of the origin, it becomes a rotation. 
The first such linearization result is
due to Siegel \cite{Sieg}, we recall it below:

\begin{thm}[\cite{Sieg}]
Suppose $\theta\in(0,1)$ is an irrational number of bounded type
\footnote{``Bounded type'' is the terminology of the renormalization theory,
the term commonly used by the number theorists is ``constant type''.}
 that is, it is 
represented by an infinite continued fraction with positive terms
$\theta=[a_1,a_2,\ldots,a_n,\ldots]$
such that $\sup a_i<\infty$. Then any analytic germ $f(z)$ with the rotation number $\theta$
is linearizable.
\end{thm}

\noindent
For an analytic mapping $f$ with a linearizable irrational indifferent fixed point, the 
maximal linearization domain is called a Siegel disk. 
In the cases when $f$ has a natural domain of definition, the structure of the boundary 
of the Siegel disk allows further study.
Let us consider the particular 
example of a quadratic polynomial 
$$f_\theta(z)=e^{2\pi i\theta}z+z^2$$ with $\theta$ an irrational
of bounded type. The following properties then hold true:

\begin{thm}[{\bf Siegel quadratics of bounded type}]
\label{bounded type}
Denote $\Delta$ the Siegel disk of $f_\theta$.
\begin{itemize}
\item[(I)] The boundary of the Siegel disk $\Delta$ is a quasicircle containing the
critical point of $f_\theta$.
\item[(II)] The Julia set of $f_\theta$ is locally connected.
\item[(III)] The Julia set of $f_\theta$ has zero area.
\end{itemize}
\end{thm}

\noindent
The second and third statements are theorems of Petersen \cite{Pet} (see also \cite{Ya1} for
a different proof).
The first statement is derived from the real {\it a priori} bounds for critical
 circle maps of
\'Swia\c\negthinspace tek \cite{Sw} and Herman \cite{Her}, via a quasiconformal surgery
procedure due to Douady, Ghys, Herman, and Shishikura. We will make use of this surgery 
further in the paper, let us therefore briefly recall it below.

\begin{proof}[Proof of (I)]
Consider the one-parameter family of degree three Blaschke products
$$Q_\tau(z)=e^{2\pi i\tau}z^2\frac{z-3}{1-3z}.$$
For each value of $\tau$ the restriction $Q_\tau|_{S^1}$ is an analytic critical
circle mapping, with critical point $c=1$ and critical value $Q_\tau(1)=e^{2\pi i\tau}$.
For each irrational number $\theta\in\TT\setminus\QQ$ let us denote 
$\tau(\theta)\in\TT$ the unique value of the parameter for which the rotation
number of $Q_{\tau(\theta)}|_{S^1}$ is equal to $\theta$ and set
\begin{equation}
\label{blaschke}
F_\theta\equiv Q_{\tau(\theta)}.
\end{equation}
Consider the rigid rotation $r_\theta(z)=e^{2\pi i\theta}z$ and let $\phi$ be the 
conjugacy
$$\phi\circ F_\theta|_{S^1}\circ\phi=r_\theta|_{S^1}$$
for which $\phi(1)=1$. By \cite{Sw} and \cite{Her} if $\theta$ is of bounded type, then
$\phi$ is a quasisymmetric map.
Let us select a  quasiconformal extension of $\phi$ to the unit disk $\DD$.
A new dynamical system 
$\tl F_\theta$ will be defined as
\begin{equation}
\label{eq surgery}
\tl F_\theta=\left\{\begin{array}{l}
                    F_\theta(z),\text{ when }z\notin\DD,\\
                    \phi^{-1}\circ r_\theta   \circ\phi(z),\text{ when }z\in\DD
                \end{array}\right.
\end{equation}
We  define a new complex structure $\mu$ in $\DD$ setting it equal to the pull-back of the
standard complex structure $\phi^*\sigma_0$. We extend $\mu$ to the outside of 
$\DD$ setting it equal to $(\tl F_\theta^n)^*\circ \phi^* \sigma_0$ for a point $z$
such that $\tl F_\theta^n(z)\in\DD$ and equal to $\sigma_0$ elsewhere.
By construction, $\tl F_\theta^*\mu=\mu$. 
By the Measurable Riemann Mapping Theorem, there exists a quasiconformal mapping
$\psi:\hat C\to\hat C$ such that
\begin{equation}
\label{qc map}
\psi^*\sigma_0=\mu \text{ a.e., and }\psi(\infty)=\infty,\;\psi(0)=0,\;\psi(1)=-e^{2\pi i\theta}/2
\end{equation}
Then
$$f_\theta=\psi^{-1}\circ\tl F_\theta\circ\psi$$
and the proof is complete.
\end{proof}

\medskip

\noindent
{\bf McMullen's results on renormalization of Siegel disks of bounded type.}
For the polynomial $f_\theta$ as in \thmref{bounded type} the restriction to $\partial \Delta$
is topologically conjugate to $R_\theta$. If $c\in\partial\Delta$ is 
the critical point of $f_\theta$, it is natural again (cf. \cite{MN,Wi,MP}) 
to define the $n$-th pre-renormalization
of $f_\theta|_{\partial \Delta}$
 as the first return map of the union of the arcs $A_n=[c,f^{q_n}_\theta(c)]$,
$A_{n+1}=[f^{q_n}_\theta(c),c]$, which is 
\begin{equation}
\label{eq1}
(f^{q_{n+1}}|_{A_n},f^{q_n}|_{A_{n+1}}).
\end{equation}
McMullen \cite{McM3} defines renormalizations of Siegel quadratics using complexified
versions of (\ref{eq1}) similar to de~Faria's holomorphic commuting pairs (cf. \cite{dF1,dF2}):

\begin{defn}[{\bf McMullen's holomoprhic pairs}] 
\label{defn1}
Let $U_1$, $U_2$, $V$ be quasidisks
in $\CC$ with $U_i\subset V$. A holomorphic pair consists of two  homeomorphisms
$g_i:\overline{U_i}\mapsto\overline{V}$ univalent on the interior, such that:
\begin{itemize}
\item $V\setminus \overline{U_1\cup U_2}$ is also a quasidisk;
\item $ \overline{U_i}\cap \partial V=I_i$ is an arc;
\item $g_i(I_i)\subset I_1\cup I_2$ for $i=1,2$;
\item $\overline{U_1}\cap\overline{U_2}=\{c\}$ is a single point.
\end{itemize}
\end{defn}
\noindent
The filled Julia set of a holomorphic pair $(g_1,g_2)$
is by definition the set $K(g_1,g_2)$ of non-escaping points for this dynamical system.

\realfig{renormfig}{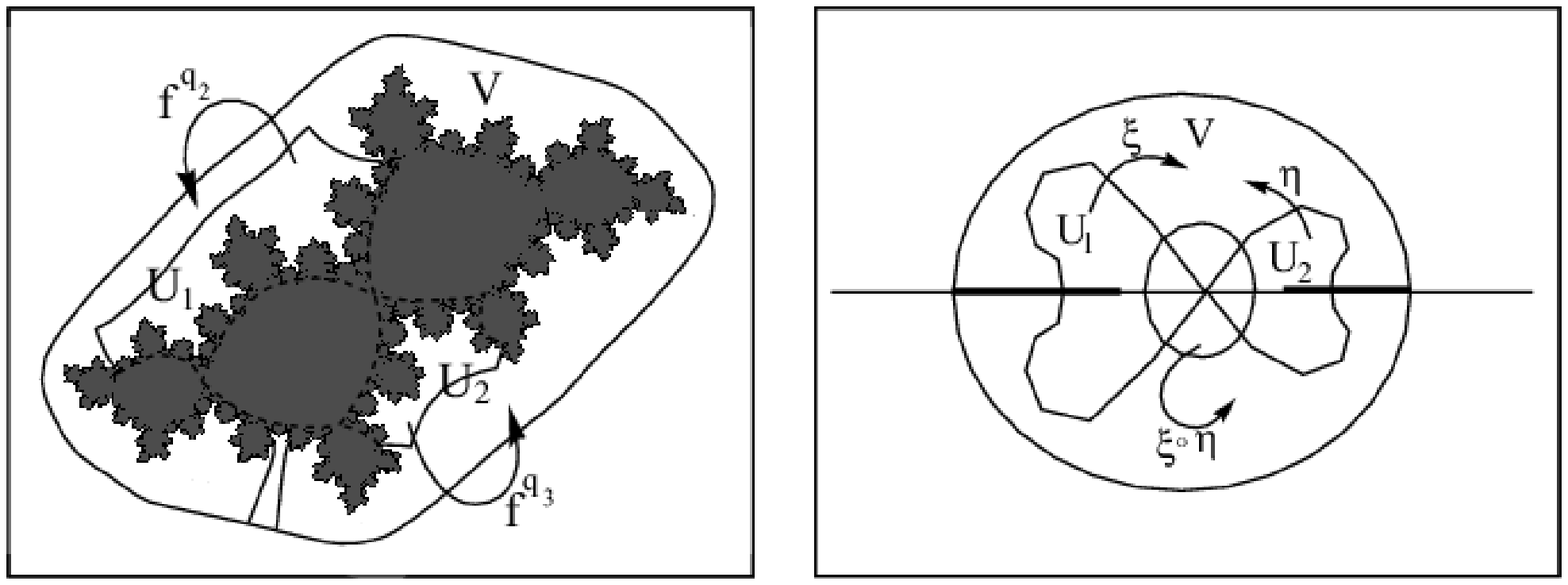}{On the left, the second pre-renormalization of the 
golden mean Siegel disk; on the right a de~Faria's holomorphic commuting pair.}{1.1\hsize}

\noindent
A pre-renormalization of $f_{\theta}$
in the sense of \cite{McM3} is a holomorphic  pair
\begin{equation}
\label{hol pair}
p\cR^n_{\text{McM}}f=(f^{q_n}:U_1^n\to V^n, f^{q_{n+1}}:U_2^n\to V^n)
\end{equation}
such that $I^n_1=A_n$ and $I^n_2=A_{n+1}$. The corresponding renormalization $\cR^n_{\text{McM}}f$ is the above
holomorphic pair linearly rescaled by a mapping sending $c$ to $0$ and the other endpoint of $A_n$
to $1$.

\begin{prop}[\cite{McM3}]
\label{mcm prepren}
A pre-renormalization (\ref{hol pair}) exists for every value of $n\in\NN$. 
\end{prop}

\noindent
This evidently follows from \thmref{complex bounds} and the surgery construction of 
\thmref{bounded type}. Below we briefly outline a different argument due to McMullen,
making use of the local connectivity
of the Julia set $J(f_\theta)$.

The domain $V^n$ is taken to be the annulus 
between the boundary of $\Delta$ and some equipotential $E$ of the Julia set of $f_{\theta}$,
from which a slight thickening of the external ray landing at $f_{\theta}(c)$ is removed.
The domain $\overline{U_1^n}$ is then the univalent pull-backs of 
$\overline{V^n}$ along the orbit
$$I_i\mapsto f(I_1^n)\mapsto\cdots\mapsto f^{q_n}(I_1^n),$$ 
and similarly for $U_2^n$. 

A holomorphic pair $\cP$ (\ref{hol pair}) under the the inverse of the surgery map
$\psi$ (\ref{qc map}) and after Schwarz reflection gives rise to a holomorphic
commuting pair $\cH=\cR^k F_\theta$. Let us say that $\cP$ has a complex bound $\mu$
if $\cH$ possesses this bound.

Let us now fix $N\in \NN$ and set 
$$\theta_*=\cfrac{1}{N+\cfrac{1}{N+\cfrac{1}{\cdots}}}$$  Denote
$p_n/q_n$ the $n$-th convergent of this infinite continued fraction.

\begin{prop}[\cite{McM3}]
\label{conj}
There exists a constant $K>1$ such that the following holds.
Every two pre-renormalizations (\ref{hol pair}) of
$f_{\theta_*}$ 
constructed as above
are $K$-quasiconformally conjugate. The conjugacy 
extends (anti)conformally to a neighborhood of $c$ in $\Delta$.
\end{prop}

\begin{thm}[\cite{McM3}]
\label{mcmullens thm}
Suppose $g$ is any analytic mapping with a golden-mean Siegel disk whose 
boundary is a quasicircle containing a single simple critical point, and let $(g_1,g_2)$ 
be a pre-renormalization of $g$ in the above sense. Assume that for some $n\in\NN$
the holomorphic pair $(g_1,g_2)$ is quasiconformally conjugate to the $n$-th
pre-renormalization of $f_{\theta_*}$ and the conjugacy extends (anti)conformally to 
a neighborhood of $c$ in $\Delta$.
Then the extended conjugacy $\psi$ is $C^{1+\alpha}$-(anti)conformal at the critical point $c$.
\end{thm}

\noindent
The proof is based on the following concept developed in \cite{McM2}:
\begin{defn}
A point $x$ is a measurable deep point of a compact $\Lambda\subset \CC$ if
there exists $\delta>0$ such that for all $r>0$
$$\text{area}(D_r(x)\setminus\Lambda)=O(r^{2+\delta})$$
\end{defn}

\noindent
The following theorem about measurable deep points appeared in \cite{McM2}:
\begin{thm}[\bf Boundary conformality]
\label{bound conf}
Let $\phi:\CC\to\CC$ be a quasiconformal map with $\bar\partial\phi=0$ on a measurable set
$\Omega$, and let $x$ be a measurable deep point of $\Omega$. Then $\phi$ 
is $C^{1+\alpha}$-conformal at $x$.
\end{thm}

\noindent
For the holomorphic pair $$g_n=(g_n^1,g_n^2):\cup U_n^i\to U_n$$ 
which is the $n$-th pre-renormalization
of $f_{\theta_*}$ in the sense of McMullen the 
filled Julia set $K(g_n^1,g_n^2)$ is a subset of $J(f_\theta)$
and hence has zero area, and no measurable deep points. McMullen considers its thickening
$$K_\eps(g_n^1,g_n^2)=K(g_n^1,g_n^2)
\cup_{\ell>0}(g_n)^{-\ell}(D_\eps(c)\setminus U_n),$$
and shows that for every $\eps>0$ the critical point $c$ is a measurable deep point of this
set. The conjugacy constructed in \propref{conj} is then $C^{1+\alpha}$-conformal
by \thmref{bound conf}.

\noindent
Denote $\cR_{\text{McM}}^n f_{\theta_*}$ the $n$-th pre-renormalization $g_n$
 rescaled linearly so that $I^n_1$ is bounded by $0$ and $1$, with $c$ mapping
to $0$.
Applying \thmref{mcmullens thm} to the particular case of these renormalized pairs, we have 

\begin{cor}
\label{cor mcm}
The renormalizations $\cR_{\text{McM}}^n f_{\theta_*}$ converge geometrically fast in $n$
with respect to the uniform metric on compact sets.
\end{cor}

\noindent
Denoting the limiting pair $(\hat \omega,\hat\upsilon)$ we see that it
 is fixed under the operation $\cR_{\text{McM}}$. Moreover, if $g$ is as in 
\thmref{mcmullens thm} then 
\begin{equation}
\label{mcm fixed pt}
\cR_{\text{McM}}^n g\to (\hat\omega,\hat\upsilon),
\end{equation}
again at a geometric rate.

\medskip

\noindent
{\bf Cylinder renormalization.}
The cylinder renormalization was introduced in \cite{Ya3}, for a detailed discussion
we refer the reader to that paper. 
Firstly, let us define some function spaces.
Denote $\pi$ the natural projection $\CC\to\CC/\ZZ$, and
 $p:\CC/\ZZ\to \bar\CC$ the conformal isomorphism given by $p(z)=e^{2\pi i z}$.
For a topological disk $W\subset\CC$ containing $0$ and $1$ we will denote 
${\aaa A}_W$ the Banach space of bounded analytic functions in $W$ equipped with the
sup norm. Let us denote ${\aaa C}_W$ the Banach subspace of 
${\aaa A}_W$ consisting of analytic mappings $h:W\to \CC$ such that
$h(0)=0$ and $h'(1)=0$. 


The cylinder renormalization operator is defined as follows. 
Let $f\in{\aaa C}_W$.
 Suppose that for $n\in\NN$ there exists a simple arc $l$ which connects a fixed point 
$a$ of $f^n$ to $0$, and has the property that $f^n(l)$ is again a simple arc 
whose only intersection with $l$ is at the two endpoints. Let $C_f$ be the topological 
disk in $\CC\setminus\{0\}$ bounded by $l$ and $f^n(l)$.
We say that $C_f$ is a {\it fundamental crescent} if 
the iterate $f^{-n}|_{C_f}$ mapping $f^{n}(l)$ to $l$  is defined and univalent, and the quotient
of $\overline{C_f\cup f^{-n}(C_f)}\setminus \{0,a\}$ by the iterate $f^n$ is conformally
isomorphic to $\CC/\ZZ$. 
Let us denote $R_f$ the first return map of $C_f$,
and let us assume that this map has a critical point $z$ corresponding to the orbit of
$0$. 
 Let $g$ be the map $R_f$ becomes under the above isomorphism, mapping $z$ to $0$,
and $h=p^{-1}\circ g\circ p$.
We say that $f$ is {\it cylinder renormalizable}, if $h\in {\aaa C}_V$ for 
some $V$, and call $h$ a {\it cylinder renormalization} of $f$.


\begin{prop}
\label{properties cyl ren}
Suppose $f\in{\aaa C}_W$ is cylinder renormalizable, and its renormalization $h_f$ is contained
in ${\aaa C}_V$. Denote $C_f$ the fundamental crescent corresponding to the
renormalization. Then the following holds.
\begin{itemize}
\item Every other fundamental crescent $C'_f$ with the same endpoints as $C_f$, and such that
$C'_f\cup C_f$ is a topological disk, produces the same renormalized map $h_f$. 
\item There exists an open neighborhood $U(f)\subset{\aaa C}_W$ such that
every map $g\in{\aaa C}_W$ is cylinder renormalizable, with a fundamental crescent $C_g$ 
which can be chosen to move continuously with $g$.
\item Moreover,  the dependence $g\mapsto h_g$ of the
cylinder renormalization on the map $g$ is an analytic mapping ${\aaa C}_W\to{\aaa C}_V$.
\end{itemize}
\end{prop}
\begin{pf}
The arguments from \cite{Ya3} apply {\it mutatis mutandis}.
\end{pf}

\noindent
Let us now concentrate on the case of Siegel quadratics with a golden rotation number.

\begin{thm}
\label{thm cyl ren siegel}
Let $f_{\theta_*}$ be as above. There exists a sequence $g_n$, $n\in\NN$ of cylinder 
renormalizations of $f_{\theta_*}$ with the following properties.
\begin{itemize}
\item[(I)] There exists an increasing sequence of natural numbers $k_n$ such that $g_n$ is a cylinder
renormalization of $f_{\theta_*}$ with period $q_{k_n}$.
For every $n$, the map
 $g_n$ has a Siegel disk with rotation number $\theta_*$ centered at the origin, 
whose boundary is a quasicircle, containing the critical
point $1$. 
\item[(II)] Let 
$$P_n:U^n_1\cup U^n_2\to V^n$$ be a sequence of holomorphic pairs (\ref{hol pair})
with a uniform complex  bound $\mu$. 
There exists $K\in\NN$ such that for every $n$ the fundamental 
crescent $C_n$ corresponding to $g_n$ is contained in the union
of the closures of the domains $U_1^{n-K}$ and $\Delta$.
\item[(III)] Finally, for $n_2>n_1$, the map $g_{n_2}$ is a cylinder 
renormalization of $g_{n_1}$.
\end{itemize}
\end{thm}

\realfig{cresc}{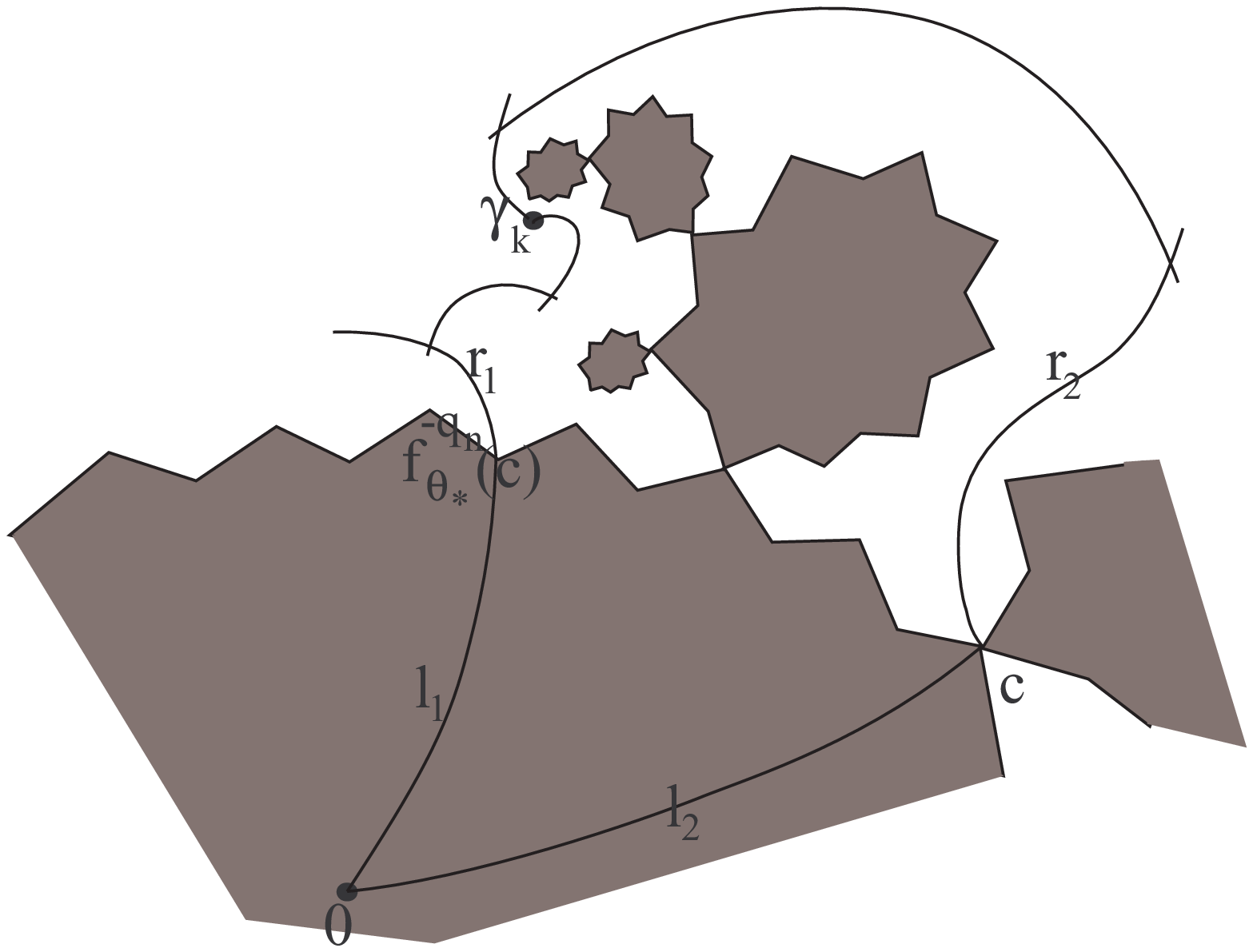}{}{10cm}

We will give two different proofs of the above theorem: 

\begin{proof}[Proof of \thmref{thm cyl ren siegel} using the surgery on a Blaschke product]
Let $F\equiv F_{\theta_*}$ be the Blaschke product (\ref{blaschke}). 
Denote $\zeta_n$ the $n$-th renormalization $\cR^n F$ of this critical circle mapping.
Each $\zeta_n$ is a commuting pair in the Epstein class, and therefore is cylinder renormalizable
(Lemma 7.6 of \cite{Ya3}). By complex {\it a priori} bounds (\thmref{complex bounds})
we can choose a fundamental crescent $W_n$ for $\zeta_n$ so that $W_{n+M}\subset W_n$ for some
fixed $M\in\NN$. We may, moreover, select $W_n$ compactly contained in the domain of the
holomorphic commuting pair $\cH_n$ of \thmref{complex bounds}.
Denote $W^+_n=W_n\cap\HH$ the upper half of the fundamental crescent. Let $p$ and $F^{q_n}(p)$ 
be the points of intersection of the boundary of $W_n$ with the unit circle.
Let $\phi:\DD\to\DD$ be the quasiconformal mapping of (\ref{eq surgery}) and $\psi$ be the change of coordinates 
(\ref{qc map}). Let $a=\phi(p)$ so that $b=r_{\theta_*}^{q_n}(a)=\phi(F^{q_n}(p))$, and denote $L\subset \DD$ the 
domain bounded by the line segments $[0,a],$ $[0,b]$ and the arc $[a,b]\subset S^1$.
Set $H_n=W_n^+\cup \phi^{-1}(L)$
and  $C_n=\psi^{-1}(H_n)$. By construction, $C_n$ is a fundamental crescent for 
$f_{\theta_*}$. Denote $h_n$ the corresponding cylinder renormalization. Setting
$g_n=h_{Mn}$ we have the required properties.

\end{proof}

\begin{proof}[Proof of \thmref{thm cyl ren siegel} using an idea of X. Buff (private communication)]
We will explicitly construct a sequence 
of fundamental crescents $C_k$ for the map $f_{\theta_*}$. 
Denote $l_1,l_2\subset\bar\Delta$ the internal rays of the Siegel disk terminating at the
endpoints of $[c,f^{-q_n}_{\theta_*}(c)]\subset\partial\Delta$,
 and let $r_1$, $r_2$ be segments of the two
 external rays of $K(f_{\theta_*})$ landing
at the same two points. 
Let $\gamma_k$ be the repelling
fixed point of $f_{\theta_*}^{q_k}$ in $U_1^k$. Let  $Q_k\ni\gamma_k$ be a  linearizing neighborhood of
$\gamma_k$, and $w_k:Q_k\to \DD$ be the linearizing coordinate, conjugating
$f_{\theta_*}^{q_k}$ to $z\mapsto \lambda_kz$ where $\lambda_k=(f_{\theta_*}^{q_k})'(\gamma_k)$.
Denote $$\Psi:\{\operatorname{Re}(z)<0\}\to\DD$$ the exponential map $z\mapsto e^{z}$,
and let $L$ be the ray 
$$\{\operatorname{Re}(z)=\operatorname{Im}(z), \;\operatorname{Re}(z)<0\}.$$
Denote $s_2$ an arc of the logarythmic spiral $w_k^{-1}\circ \Psi(L)$ which terminates
at $\gamma_k$ and whose other endpoint lies outside of the filled Julia set $K(f_{\theta_*})$.
Set $s_1$ to be the component of the preimage $f_{\theta_*}^{-q_k}(s_2)$ which contains
$\gamma_k$.
Let $t_2$ be a curve consisting of an arc of an external ray and of an equipotential
connecting $s_2$ with $r_2$, and $t_1$ its $f^{-q_k}_{\theta_*}$-preimage
connecting $s_1$ with $r_1$. 
By construction,
the curves $L_i=l_i\cup r_i \cup t_i\cup s_i$ are disjoint and $f^{q_k}_{\theta_*}(L_1)=L_2$.
The local pictures of dynamics at $0$ and at $\gamma_k$ imply that 
 $L_1$ and $L_2$ bound a fundamental crescent $C_k$ of $f_{\theta_*}$
(cf. \figref{cresc}). The map $h_k$
which is the corresponding cylinder renormalization satisfies the property (I).
By \corref{cor mcm} there exists $M\in\NN$ such that for every $k$ the crescent
$C_{M+k}\subset C_k$. Setting $g_k=h_{Mk}$ we have the desired properties.
\end{proof}

\subsection{Parabolic renormalization}
For further reference,
we give here a brief definition of the relevant version of parabolic renormalization.
We refer the reader to \cite{Sh}, and also to \cite{EY}. The latter reference contains a detailed
discussion of parabolic renormalization for critical circle mappings.

Let $f\in{\aaa C}_U$ have a fixed point $p$ with eigenvalue $1$.
Assume further, that $p$ is  a simple parabolic point, and denote
$C^A$ and $C^R$ the attracting and repelling
Fatou cylinders of $p$. The overlap of the attracting and repelling petals
of $f$ induces an analytic {\it {\'E}calle-Voronin mapping} $\cE_f$ from open neighborhoods
$W^+$, $W^-$ of the ends $\oplus,\;\ominus$ of $C^R$ to $C^A$. We normalize the situation by
requiring that $(\cE_f)'(\oplus)=1$. 

An arbitrary choice of an affine isomorphism $\tau:C^A\to C^R$ induces a dynamical
system
$$g_\tau\equiv \tau\circ\cE_f:W^+\cup W^-\to C^R.$$
Denote $e(z)=e^{2\pi i z}$ the conformal isomorphism $C^A\to \CC\setminus\{0\}$,
sending the end $\oplus$ to the puncure at the origin.
Now let $[a_1,a_2,\ldots,a_n,\ldots]$ be any formal continued fraction with 
$a_i\in\NN\cup \{\infty\}$. If there is no symbol $\infty$
present in the sequence $[a_i]$, then the continued fraction 
converges to a well-defined $\theta\in\RR\setminus\QQ$. Otherwise, 
assuming that $i$ is the first position in which $\infty$ is encountered,
we let $\theta$ to be the rational $[a_1,\ldots,a_{i-1}]$.

We set
$$\cP_\theta f=e\circ g_{\tau(\theta)}\circ e^{-1},$$ 
for the unique choice of $\tau(\theta)$ such that $(\cP_\theta f)'(0)=e^{2\pi i\theta}$, and
call it the {\it parabolic renormalization} of $f$ corresponding to $\theta$.

Parabolic renormalization can be seen as the limiting case of cylinder
renormalization as seen from the following (cf. \cite{Sh,EY}):

\begin{thm}
\label{continuity fatou}
For $f$ as above, there exists a neighborhood $\cU(f)\subset {\aaa C}_U$ such that the
following holds.
Every $h\in\cU(f)$ with $h'(0)=e^{2\pi it}$, $t\neq 0$ is cylinder renormalizable with
period $1$.

Moreover, fix $\theta=[a_1,a_2,\ldots,a_n,\ldots]$ with
$a_i\in\NN\cup \{\infty\}$. Let $h_i\in \cU(f)$ converge to $f$, and have $(h_i)'(0)=e^{2\pi i t_i}$.
When $t_i\neq 0$, let $g_i$ be the cylinder renormalization of $h_i$ of period $1$. When 
$t_i=0$, set $g_i=\cP_\theta h_i$. 

Assume that $(g_i)'(0)\to e^{2\pi i\theta}$. Then $$g_i\to \cP_\theta f$$
uniformly in some neighborhood of $0$.
\end{thm}

\section{Proof of the main theorem}

\noindent
{\bf The definition of $\cren$ and the existence of a fixed point $\hat f$.}
\begin{prop}
\label{pr1}
Let $g_n$ be the sequence of cylinder renormalizations of $f_{\theta_*}$ constructed
in \thmref{thm cyl ren siegel}. 
There exists $M\in\NN$ and a domain $U\supset\{0,1\}$ such that for $n\geq M$ the renormalizations $g_n\in {\aaa C}_U$, 
and converge geometrically fast in the uniform topology
 to a map $\hat f\in{\aaa C}_U$. The map $\hat f$ has a Siegel disk 
$\Delta_{\hat f}$ which is compactly contained in $U$.
\end{prop}
\begin{pf}
Let  $K$ be as in \thmref{thm cyl ren siegel}.
By \thmref{mcmullens thm} $\cR^{Kn}_{\text{McM}}f$ converge geometrically fast to a holomorphic
pair $(\hat \omega,\hat\upsilon)$ which is fixed under $\cR_{\text{McM}}$.
This and \thmref{thm cyl ren siegel} imply the desired claim.
\end{pf}

\begin{prop}
\label{pr2}
There exists a topological disk $W\Supset \Delta_{\hat f}$ such that the following holds.
For every topological disk $V$ with $W\Supset V\Supset\Delta_{\hat f}$ we have:
\begin{itemize}
\item[(I)] 
Denote $g$ the restriction $\hat f|_V$. Then the map $g$ has a cylinder renormalization $\tl g$ 
whose domain of definition compactly contains $W$ and such that $\tl g|_W\equiv \hat f|_W$.
\item[(II)] Moreover, there exists $n$ such that the $n$-th pre-renormalization in the
sense of McMullen (\ref{hol pair})
$$p\cR^n_{\text{McM}}:U_1^n\cup U_2^n\to V^n $$
 exists for $g$.
\item[(III)] Finally, let $J_i$ be an element of the $n$-th dynamical partition (\ref{dyn part})
of $g|_{\partial\Delta_{\hat f}}$ corresponding to the orbit of $0$,
and $g^k(J_i)\subset \overline{U_1^n\cup U_2^n}$.
Then the inverse branch $\tau_i$ mapping $g^k(J_i)$ to $J_i$
 univalently extends to $(V^n\setminus \partial\Delta_{\hat f})\cup g^{k}(J_i)$
and maps it to a subset of $V$.
\end{itemize} 
\end{prop}
\begin{pf}
We will make use of the quasiconformal conjugacy $\psi=\psi_{\theta_*}$ (\ref{qc map}) 
conjugating the quadratic map $f_{\theta_*}$ to the modified Blaschke product $\tl F_{\theta_*}$.
Let $V$ be any domain compactly containing $\Delta_{\hat f}$. Set $X=\psi(V)\setminus \overline{\DD}$, and 
$Y$ to be the union of $X$ with its reflection in $S^1$, together with the arc $\psi(V)\cap S^1$.
 Let $G=F_{\theta_*}|_Y$.
By complex {\it a priori} bounds (\thmref{complex bounds})
 there exists $N$ such that for every analytic critical circle map $h$
in $Y$ with a single critical point,  for $n\geq N$ the commuting pair
$ p\cR^n G$ extends to a holomorphic
commuting pair $\cH_h^n\in\hol(\mu)$. By a compactness argument, $N=N(Y)$.

The pull-back of the ``upper half'' of the holomorphic commuting pair $\cH_G^n$ by the
conjugacy $\psi$
is a holomorphic pair 
$$p\cR^n_{\text{McM}}:U_1^n\cup U_2^n\to V^n,$$
which implies part (II). 
Note that the universality of the bound $\mu$ 
and \thmref{thm cyl ren siegel} (II) imply that there exists $W$ independent of the initial
choice of $V$ such that (I) holds.
 Part (III) follows
from the same considerations.

\end{pf}

\noindent
\begin{defn}
Let us select a subdomain $V\Subset W$, and $n\in\NN$ as in the previous Proposition. Let $U\Subset W$ be an
open neighborhood of $\Delta_{\hat f}$ such that 
\begin{equation}
\label{cr1}
U\Subset \cup\tau_i(V^n)\cup\overline{\Delta_{\hat f}}\text{, and }
U\Supset \cup\tau_i(U_i^n)
\end{equation}
By \propref{pr2} (I) there exists a  cylinder renormalization transforming
$\hat f|_U$ to $\hat f|_W$. By \propref{properties cyl ren} it extends to an analytic operator
from an open subset of ${\aaa C}_U$ to ${\aaa C}_U$. We will denote this operator $\cren$,
and call it the cylinder renormalization operator.
\end{defn}

\medskip
\noindent
{\bf The expanding direction of} $\cren$.
Let us denote $\cL=D_{\hat f}\cren:{\aaa C}_U\to {\aaa C}_U$.
We first establish that:

\begin{prop}
The operator $\cL$ is compact.
\end{prop}
\begin{pf}
Denote $B_1$ the unit ball in ${\aaa C}_U$ and let $v\in B_1$. 
By definition of $\cren$, 
the vector field
$\cL v$ is an analytic vector field in the  domain $W\Supset U$. 
Denote $C\Subset U$ a fundamental crescent which corresponds to the cylinder renormalization
operator, and let $\Phi:C\to \hat \CC$ be its uniformization.
 The first return map $R_C$ of $C$ under $\hat f|_U$ is a bounded piecewise
analytic map. The restriction of $R_C$ to $\Phi^{-1}(U)$ is a finite collection
of iterates $\hat f^j$, and the compactness considerations imply that there
exists $C$ independent
of $v$ such that $||\cL v|_U||<C$. Since the $C$-bounded functions in ${\aaa C}_U$
which analytically extend to $W$ form a compact set, the image
$\overline{\cL({B_1})}$ is compact.
\end{pf}

\begin{prop}
The operator $\cL$ has an eigenvalue $\lambda$ with $|\lambda|>1$.
\end{prop}
\begin{pf}
Let $v(z)$ be a vector field in ${\aaa C}_U$,
$$v(z)=v'(0)z+o(z).$$
Denote $\gamma_v$ the quantity 
$$\gamma_v=\frac{v'(0)}{\hat f'(0)}={e^{-2\pi i\theta_*}}v'(0)$$
For a smooth family
$$\hat f_t(z)=\hat f(z)+tv(z)+o(t),$$
we have 
$$\hat f_t(z)=\alpha^v_t(z)(\hat f'(0)z+o(z)),\text{ where }\alpha_t(0)=1+t\gamma_v+o(t).$$
The $q_{m+1}$-st iterate
$$\hat f_t^{q_{m+1}}(z)=(\alpha^v_t(z))^{q_{m+1}}((\hat f'(0))^{q_{m+1}}z+o(z)).$$
In the neighborhood of $0$ the renormalized vector field $\cL v$ is obtained
by applying a uniformizing coordinate 
$$\Psi(z)=(z+o(z))^\beta,\text{ where }\beta=\frac{1}{\theta_* q_m\mod 1}.$$
Hence,
$$\alpha_t^{\cL v}(0)=[(\alpha_t^v(0))^{q_{m+1}}]^\beta,$$
so
$$\gamma_{\cL v}=\Lambda \gamma_v,\text{ where }\Lambda=\beta q_{m+1}>1.$$
Hence the spectral radius
$$R_{\text{Sp}}(\cL v)>1,$$
and since every 
 non-zero element of the spectrum of a compact operator is
an eigenvalue, the claim follows.
\end{pf}

\medskip
\noindent
{\bf Stable direction.}
Denote $\cW\subset{\aaa C}_U$ the collection of mappings $f$ such that
$f'(0)=e^{2\pi i \theta_*}$.

\medskip
\noindent
{\bf Inou and Shishikura's results.} Inou and Shishikura \cite{Shi} have recently established the following result:
\begin{thm}
\label{thm-shi}
There exist $N_0\in \NN$, a pair of topological disks $\widetilde W\Supset W\ni 0$, and an open neighborhood
$\cV\subset {\aaa C}_W$ so that the following is true.

\begin{itemize}
\item Let $\theta=[a_1,a_2,\ldots]$ with $a_i\in \NN\cup\infty$ and $a_i\geq N_0$. For every $f\in \cV$ with
$f'(0)=e^{2\pi i\theta}$ we have the following. 
If $a_1\neq \infty$, then $f$ is cylinder renormalizable
with period $1$, and the corresponding cylinder renormalization $g\in \cV$.
Otherwise, $$\cP_{[a_2,a_3,\ldots]}f\in \cV\cap {\aaa C}_{\widetilde W}.$$

\item Moreover, consider the quadratic polynomial $f=f_\theta(z)$. Set $g_i$ to 
be the sequence of cylinder/parabolic renormalizations of $f$ with 
$$(g_i)'(0)=e^{2\pi i[a_{i+1},a_{i+2},\ldots]}.$$
Then there exists $j\in\NN$ such that $g_j\in\cV$.

\item Finally, the neighborhood $\cV$ can be chosen sufficiently small, so that, in particular,
 for every $f\in \cV$ the critical point $1$ is not fixed.

\end{itemize}
\end{thm}

\noindent
Note that the claim (V) of the Main Theorem immediately follows from the above result.

\begin{prop}
\label{prop3}
There exists an open neighborhood $\cU\subset{\aaa C}_U$ of $\hat f$ such that
there exists $N\in\NN$ such that
for every $f\in\cW\cap \cU$ the following holds:
\begin{itemize}
\item[(I)] the boundary of the Siegel disk $\partial\Delta_f$ is a quasicircle containing the
critical point $1$;
\item[(II)]
 the pair of mappings $$(f^{q_N}|_{[1,f^{q_{N+1}}(1)]},f^{q_{N+1}}|_{[1,f^{q_{N}}(1)]})$$
extends to a holomorphic pair, which is quasiconformally conjugate to $(\hat\omega,\hat\upsilon)$ (\ref{mcm fixed pt})
in such a way, that the conjugacy conformally extends to the Siegel disks.
\end{itemize}
\end{prop}
\begin{pf}
From \thmref{thm-shi} there exists a neighborhood $\cU$ such that for every $f\in \cU$ the
critical orbit $\{f^n(1)\}$ is infinite and compactly contained in $\cU$. Since the
dependence $f\mapsto f^n(1)$ is holomorphic,  the Bers-Royden
Theorem implies that the closure $\overline{ \{f^n(1)\} }$ is a $K$-quasicircle $\Upsilon_f$ for some $K=K(\cU)>1$.
Given the invariance of $\Upsilon_f$, and considerations of the Denjoy-Wolff Theorem, we see that
$\Upsilon_f=\partial \Delta_f$.
The standard pull-back argument implies that $f$ is $K_1=K_1(\cU)$-quasiconformally conjugate to 
$\hat f$ on a subdomain $V\subset U$ with $\Delta_{ f}\Subset V$.

By \propref{pr2} the neighborhood $\cU$ may be chosen small enough so that
(II) holds for some $N\in\NN$.
\end{pf}

\begin{prop}
There exists a codimension one subspace $T$ of ${\aaa C}_U$ such that 
$$R_{\text{sp}}(\cL|_T)<1.$$
\end{prop}

\begin{pf}
We take $T$ to be the tangent space to the codimension one submanifold
 $\cU\cap\cW\subset{\aaa C}_U$
 at the point $\hat f$. By McMullen's \thmref{mcmullens thm}, for every map $f\in \cU$,
$$\cren^k f\to\hat f$$
geometrically fast.
The Spectral Theorem for compact operators implies that
$$R_{\text{sp}}(\cL|_T)<1.$$
\end{pf}


\appendix

\section{Application to the Measure Problem}

The question of existence of polynomial Julia sets with positive area
has long been a central problem in Holomorphic Dynamics. Some years ago
Douady has outlined a program for constructing such sets in the family
of quadratic polynomials. Ch{\'e}ritat \cite{Ch} in his thesis made
a major progress in the program, reducing it to several renormalization-related 
conjectures. Very recently,
Buff and Ch{\'e}ritat \cite{BC} have announced the construction of some
positive measure examples.  In their approach, they used the work of Inou \& Shishikura \cite{Shi}
to resolve the conjectures of Ch{\'e}ritat.

In this Appendix, we demonstrate how the same result follows from hyperbolicity
of Siegel renormalization (\thmref{main thm}). Given that part (VI) of \thmref{main thm}
also relies on the work of Inou and Shishikura, this approach is not much different
from that taken by Buff and Ch{\'e}ritat. The use of the renormalization theorem, however,
makes for a shorter proof, and clarifies the connection between hyperbolicity of $\cren$
and positive measure. The latter could be  particularly useful in the further study of
measure of Cremer and Siegel Julia sets by renormalization techniques.

We formulate a theorem which by \cite{Ch} implies that some quadratic Julia sets  have
positive measure:

\begin{thm}
\label{meas-thm}
Let $[a_1,\ldots,a_k]$ be any finite continued fraction with $a_i\geq N_0$ from \thmref{thm-shi}.
Fix $N\geq N_0$.
Denote $\theta_j$ the finite continued fraction 
$$\theta_j=[a_1,\ldots,a_k,\underbrace{N,N,\ldots,N}_j,\infty].$$
Let $\CC/\ZZ\simeq C^R$ be the repelling cylinder of the quadratic polynomial $f_{\theta_j}$.
Denote $\hat K_j\subset C^R$ the projection of $K(f_{\theta_j})$, and let $C_j$ be
its complement $C^R\setminus \hat K_j$. 
Then the following holds:
\begin{itemize}
\item[(I)]  Denote 
$$h_j=\sup (\Im z_1 -\Im z_2),\text{ for } z_1,z_2\in C_j.$$
Then $\sup h_j<\infty.$ 
\item[(II)] 
$$\operatorname{Area}(C_j)\to 0.$$
\end{itemize}
\end{thm}

Note that (I) follows immediately from \thmref{thm-shi} and compactness considerations. We will 
proceed to deriving (II).

Jellouli and Ch{\'e}ritat have shown (see \cite{Ch}):
\begin{thm}
\label{sieg approx}
For every finite continued fraction  $[a_1,\ldots,a_k]$ denote 
$$\theta_j=[a_1,\ldots,a_k,\underbrace{N,N,\ldots,N}_j,\infty], \text{ and }\theta_\infty=\lim_{j\to\infty}\theta_j.$$
Then 
$$\operatorname{Area} K(f_{\theta_\infty})\setminus K(f_{\theta_j})\longrightarrow 0.$$
\end{thm}

We also recall the following fact proven in \cite{Ya1}:

\begin{lem}
\label{yam-lem}
Let $\theta_\infty$ be as above, and denote $\Delta=\Delta_{\theta_\infty}$ the Siegel disk of the quadratic
polynomial $f_{\theta_\infty}$. Then 
$$\liminf_{\eps\to 0}\frac{\operatorname{\meas} U_\eps(\Delta)\setminus K(f_{\theta_\infty})}{\operatorname{\meas} U_\eps(\Delta)}=0.$$
\end{lem}

\noindent
We remark that by later work of McMullen \cite{McM3}, in the above expression $\liminf$ can be replaced with $\lim$, and in fact, the
convergence occurs at a geometric rate. However,  this stronger result will not be needed.

The idea of the proof of \thmref{meas-thm} is rather straightforward. 
As seen from the above lemma, the basin of the Siegel disk of $f_{\theta_\infty}$ occupies most of the 
area in a sequence of arbitrarily small disks around the critical point. Hence, for the renormalization
fixed point $\hat f$, the Siegel basin has full measure in the cylinder. 
As implied by the hyperbolicity of renormalization, all sufficiently high cylinder renormalizations
of $f_{\theta_j}$ are contained in a small neighbourhood of $\hat f$. 
From this we will conclude
that the parabolic basin of $f_{\theta_j}$ on a sufficiently small scale contains most of the Siegel
basin of $\hat f$ with thin cusps removed. By the result of Jellouli and Ch{\'e}ritat the 
total area of the removed cusps can be made arbitrary small by selecting large enough $j$.
Hence, most of the measure in a fundamental crescent of $f_{\theta_j}$ will be occupied by the
parabolic basin. Now we proceed to formalize this discussion.

\begin{proof}[Proof of \thmref{meas-thm}]
First of all, observe that by \lemref{yam-lem}, and \thmref{main thm}, for every $\eps>0$
 and every disk 
$D\subset \CC/\ZZ$  there exists $K_0\in\NN$ such that for every $K>K_0$ there 
exists $M\in\NN$ for which the following holds.
Denote $g_K=\cren^K(f_{\theta_\infty})$, and let $B_K$ be the projection on the domain of $g_K$ of
$(f_{\theta_\infty})^{-M}(\Delta)$. Then 
$$\frac{\operatorname{Area}(D\setminus B_K)}{\operatorname{Area}(D)}<\eps.$$

Now select a point $z\in C_j$ and let $\zeta$ be any lift of this point. 
Let us fix $\eps>0$ and $M\in\NN$ and select a sufficiently large $j\in\NN$ and  $K+M<j$ such that 
the above inequality holds. 
Let $(\omega:U_1\to V,\upsilon:U_2\to V)$ be a McMullen renormalization of $g_K$.

\noindent
Given that the point $\zeta$ is not in the filled Julia set of $f_{\theta_j}$, 
we necessarily have one of the following two scenarios:

\begin{itemize}
\item[(*)] either there exists $n\in\NN$ such that $\zeta_n\equiv (f_{\theta_j})^n(\zeta)$ is 
contained in the domain $V\setminus (U_1\cup U_2)$;

\item[(**)] or there exists $n$ such that $\zeta_n\in\Delta$.

\end{itemize}

In the case (*), consider a disk $D_n\ni\zeta_n$ such that 
$$\frac{\operatorname{Area}(D_n\setminus f^{-M}(\Delta))}{\operatorname{Area}(D_n)}<\eps,\text{ and }
\dist(D_n,\text{Postcrit}(f_{\theta_\infty}))>\diam (D_n).$$
Denote $D\ni \zeta$ its preimage, and $\tl D\ni z$ the projection to $C^R$.
By Koebe's Distortion Theorem and \thmref{sieg approx} we have $$\operatorname{Area}(\tl D\cap C_j)/\operatorname{Area}(\tl D)<\delta(\eps)\text{ with }\delta(\eps)\underset{\eps\to 0}{\longrightarrow}0,$$
where $\operatorname{Postcrit}(f_{\theta_\infty})$ denotes the postcritical set.

In the case (**) we can select a disk $D_{n-1}\ni \zeta_{n-1}$ with the same properties, and the
same considerations apply.

The statement of the theorem follows.

\comm{By Koebe Distortion Theorem, and \thmref{sieg approx} for every $\delta>0$ and every $\eps>0$, there exists $J\in\NN$
such that for all $j>J$
we have the following.  Let $z\in C_j$ and $W\ni z$ be a Euclidean disk with diameter in the interval $[\delta,1]$.
Then
$\operatorname{Area}(W\cap C_j)<\eps$.
}

\end{proof}

\end{document}